\documentclass[letterpaper, 10 pt, conference]{ieeeconf}  % Comment this line out if you need a4paper

\IEEEoverridecommandlockouts                              % This command is only needed if 
\usepackage{graphicx} % for pdf, bitmapped graphics files
\usepackage{epstopdf}

\usepackage{subfigure}

\usepackage{caption}

\usepackage{amsmath}
\usepackage{amssymb}
\usepackage{booktabs} 

\newcommand{\barr}{\begin{array}}
	\newcommand{\earr}{\end{array}}

\newcommand{\matlab}{\textsc{Matlab}}

\newcommand{\Dc}{D_{\mathrm{c}}}
\newcommand{\Hc}{H_{\mathrm{c}}}
\newcommand{\hc}{h_{\mathrm{c}}}

\newcommand{\OR}{$\Vert$}

\usepackage{color}

\usepackage{algorithm}
\usepackage{algpseudocode}
\algrenewcommand\algorithmicrequire{\textbf{Input:}}
\algrenewcommand\algorithmicensure{\textbf{Output:}}

\newcommand{\software}[1]{\texttt{#1}}

\usepackage{float}
\newfloat{algorithm}{t}{lop}

\usepackage{color}
\definecolor{wheat}{rgb}{0.96,0.87,0.70}

\title{\LARGE \bf
A Structure Exploiting Branch-and-Bound Algorithm for Mixed-Integer Model Predictive Control
}

\author{Pedro Hespanhol$^{1}$, Rien Quirynen$^{1}$, Stefano Di Cairano$^{1}$% <-this % stops a space
	\thanks{$^{1}$Control and Dynamical Systems, Mitsubishi Electric Research Laboratories, Cambridge, MA, 02139, USA.
		{\tt\small quirynen@merl.com}}
}

\begin{document}

\maketitle
\thispagestyle{empty}
\pagestyle{empty}

%%%%%%%%%%%%%%%%%%%%%%%%%%%%%%%%%%%%%%%%%%%%%%%%%%%%%%%%%%%%%%%%%%%%%%%%%%%%%%%%
\begin{abstract}
Mixed-integer model predictive control~(MI-MPC) requires the solution of a mixed-integer quadratic program~(MIQP) at each sampling instant under strict timing constraints, where part of the state and control variables can only assume a discrete set of values. Several applications in automotive, aerospace and hybrid systems are practical examples of how such discrete-valued variables arise. We utilize the sequential nature and the problem structure of MI-MPC in order to provide a branch-and-bound algorithm that can exploit not only the block-sparse optimal control structure of the problem but that can also be warm started by propagating information from branch-and-bound trees and solution paths at previous time steps. We illustrate the computational performance of the proposed algorithm and compare against current state-of-the-art solvers for multiple MPC case studies, based on a preliminary implementation in \matlab~and C~code.
%We show how this can be done effeficiently in the context of MPC and illustrate it's potential use in embedded application by presenting the performance of the algorithm, compared to state-of-the-art commercial solvers.
\end{abstract}

%%%%%%%%%%%%%%%%%%%%%%%%%%%%%%%%%%%%%%%%%%%%%%%%%%%%%%%%%%%%%%%%%%%%%%%%%%%%%%%%
\section{INTRODUCTION}

Optimization based control and estimation techniques, such as model predictive control~(MPC) and moving horizon estimation~(MHE), allow a model-based design framework in which the system dynamics and constraints can directly be taken into account~\cite{Mayne2013}. This framework can be further extended to hybrid systems~\cite{bemporad1999control}, providing a powerful technique to model a large range of problems, e.g., including dynamical systems with mode switchings or quantized control, problems with logic rules or no-go zone constraints. 
%The first step is to use direct optimal control methods (popular examples of this include the direct multiple shooting method~\cite{bock1984multiple} and direct collocation~\cite{Betts2010,Biegler1984})  to transform the continuous time system into a (discrete time) optimization problem. 
However, the resulting optimization problems are highly non-convex because they contain variables that only take integer values. When using a quadratic objective in combination with linear system dynamics and linear inequality constraints, the resulting optimal control problem~(OCP) can be formulated as a mixed-integer quadratic program~(MIQP). 

%We formulate our discrete time direct optimal control problem as a mixed-integer quadratic program (MIQP). One can obtain such QP after utilizing, for example, an explicit fixed-step integration scheme to go from the continuous-time OCP to the MIQP formulation. 
%We let $[N]$ denote the set of indices $\{0,\ldots,N\}$ denoting the discrete time stages. 
We aim to solve MIQP problems of the following form:
\begin{subequations} \label{MI-OCP}
	\begin{alignat}{5}
	\underset{X,\,U}{\text{min}} \quad &\frac{1}{2}\sum_{i=0}^{N-1} x^{\top}_{i}Q_{i}x_i + u^{\top}_{i}R_{i}u_i + x^{\top}_{N}Px_{N}&& \label{MS:obj}\\
	\text{s.t.} \quad\; & x_0 - \hat{x}_0 \;= \;0,  	\label{MS:initial}\\
	& A_{i}x_{i} + B_{i}u_{i} + a_{i} \;=\; x_{i+1}  , \quad && \hspace{-2em} i \in \{0,\ldots,N-1\}, \label{MS:Kstep}\\
	& l^{c}_{i} \leq C_{i}x_{i} + D_{i}u_{i} \leq u^{c}_{i}, \quad &&\hspace{-2em} i \in \{0,\ldots,N-1\}, \label{MS:path} \\
%	& \bar{l}^{c}_{i} \leq C_{i}x_{i} + D_{i}u_{i} \leq \bar{u}^{c}_{i}, \quad &&i \in \{0,\ldots,N\} \label{MS:path} \\
%	& E_ix_{i} + G_{i}u_{i} \in \{\bar{l}_{i} , \bar{u}_{i}\} \quad && i \in \{0,\ldots,N-1\} \label{MS:int} 
	& F_iu_{i} \in \{0, 1\}, \quad && \hspace{-2em} i \in \{0,\ldots,N-1\}, \label{MS:int} \\
	& l^{c}_{N} \leq C_{N}x_{N} \leq u^{c}_{N}, \label{MS:term}
	\end{alignat}
\end{subequations}
%We let $[N]$ denote the index set $\{0,\ldots,N\}$.
where the optimization variables are the state $X=[x_0^\top,\ldots, x_N^\top]^\top$ and control trajectory $U=[u_0^\top,\ldots, u_{N-1}^\top]^\top$. The set of constraints~\eqref{MS:int} are binary equality constraints, since the left-hand side needs to be equal to either $0$ or $1$. For simplicity of notation, we further consider only binary control variables instead of more general integer constraints for an affine function of both state and control variables.
%We note that enforcing a subset of he variables to be binary, or  to belong in a certain range of integers are only special cases of the constraints in~\eqref{MS:int}.   
%Note that MPC for any hybrid system can be formulated as in~\eqref{MI-OCP}, by using a mixed logical dynamical~(MLD) model reformulation~\cite{bemporad1999control}. 
MPC for several classes of hybrid systems can be straightforwardly formulated as in~\eqref{MI-OCP}. Notable examples are mixed logical systems~\cite{bemporad1999control}, where auxiliary continuous and discrete variables can be added to the input vector.
Moreover, in combination with the binary constraints~\eqref{MS:int}, the affine inequalities~\eqref{MS:path} can model various complicated but practical restrictions on the feasible region, such as no-go zones and disjoint polyhedral constraints for states and inputs.

A hybrid MPC controller aims to solve the MIQP~\eqref{MI-OCP} at every sampling time instant. This is a difficult task, given that mixed-integer programming is $\mathcal{NP}$-hard in general, and several methods for solving such a sequence of MIQPs have been explored in the literature. These approaches can be divided into heuristic techniques, which seek to efficiently find sub-optimal solutions to the problem, and optimization algorithms which attempt to solve the MIQPs to optimality. Examples of the former include rounding and pumping schemes~\cite{achterberg2007improving,10.1007/978-3-642-29210-1_12}, approximate optimization algorithms~\cite{diamond2016general,naik2017embedded}, and approximate dynamic programming~\cite{Stellato2017}. The downside of fast heuristic approaches is often the lack of guarantees for finding an optimal or even an integer-feasible solution. Heuristic rounding-based approaches to mixed-integer nonlinear OCPs can be found, e.g., in~\cite{Kirches2010y,Sager2008}.

%that aside from not guaranteeing the an optimal solution is found, convergence to a feasible solution is not certain. 
As for solving these problems to optimality, most of the optimization algorithms for MIQPs are based on the classical branch-and-bound~(B\&B) technique~\cite{floudas1995nonlinear}. For the purpose of mixed-integer MPC, the standard B\&B strategy has been combined with various methods for solving the relaxed convex QPs. For example, a B\&B algorithm for mixed-integer MPC~(MI-MPC) has been proposed in combination with a dual active-set solver in~\cite{axehill2006mixed}, with an interior point algorithm in~\cite{frick2015embedded}, dual projected gradient methods in~\cite{Axehill17355,naik2017embedded}, a nonnegative least squares solver in~\cite{Bemporad2018}, and the alternating direction method of multipliers~(ADMM) in~\cite{stellato2018embedded}. Branch-and-bound methods for solving mixed-integer nonlinear OCPs have also been studied, e.g., in~\cite{Gerdts2005}.

Another important research topic focuses on general pre-processing and modeling techniques to reduce the size and strengthen the mixed-integer problem formulations~\cite{nemhauser1988integer}. These \emph{presolve} techniques are vital to the good performance of current state-of-the-art mixed-integer solvers~\cite{achterberg2016presolve}, such that these methods can often solve seemingly intractable problems in practice. Lastly, the branch-and-bound method itself has been extensively studied with several improvements in branching and variable selection techniques~\cite{achterberg2005branching,le2017abstract}, including recent developments in applying machine learning techniques in order to learn ``better'' branching rules~\cite{balcan2018learning}. Finally, the branch-and-bound strategy has been generalized further, e.g., using cutting planes to tighten the convex problem relaxations, resulting in \emph{branch-and-cut} or \emph{branch-and-price} variants of the algorithm~\cite{floudas1995nonlinear,nemhauser1988integer}.
Unlike state-of-the-art mixed-integer solvers, e.g., \software{GUROBI}~\cite{gurobi} and \software{MOSEK}~\cite{mosek}, our aim is to propose a tailored algorithm and its solver implementation for fast embedded MI-MPC applications, i.e., running on microprocessors with considerably less computational resources and available memory. The optimization algorithm should be relatively simple to code with a moderate use of resources, while the software implementation is preferably compact and library independent.
%\todo{I do not think we never say our goal, which is obtaining a relatively fast MI-MPC for embedded applications, i.e., library independent, simple code, compact, and moderate use of resources....
%	... this should be somewhere in the introduction (maybe also contrasting with state of the art solvers...)}

In this paper, our first contribution is to propose a branch-and-bound based MPC algorithm, which exploits the features of a recently proposed structure-exploiting primal active-set solver called \software{PRESAS}~\cite{PRESAS}. The latter algorithm is tailored to efficiently solve QPs with a block-sparse optimal control structure. 
%The performance of this solver, in comparison to other state-of-the-art QP solvers can be found in \cite{PRESAS}.
Our second contribution is to bring various mixed-integer programming techniques, such as bound strengthening, domain propagation, and advanced branching rules, to the context of MI-MPC. In particular, we present an algorithm that exploits the sequential nature of MPC, in order to warm-start the branch-and-bound search tree and to re-use information gathered at previous time steps. A similar type of approach was proposed recently by~\cite{Bemporad2018}, but in this work we provide not only a warm-start procedure for the integer variables but we also show how to improve the branching strategy by warm starting and how to efficiently combine this with presolving techniques for MI-MPC. Finally, the computational performance of the proposed algorithm, for a preliminary implementation in \matlab~and C~code, is illustrated and compared against current state-of-the-art solvers for multiple MPC case studies.

The paper is organized as follows. Section~\ref{sec:MIQP} presents the basic idea of branch-and-bound methods for mixed-integer programming. Then, Section~\ref{sec:PRESOLVE} presents presolve techniques in the context of mixed-integer optimal control. The resulting MI-MPC algorithm and its tailored warm-starting strategies are discussed in Section~\ref{sec:three} and its performance is illustrated based on multiple numerical case studies in Section~\ref{sec:caseStudies}.
Finally, Section~\ref{sec:concl} concludes the paper.
%Section~\ref{sec:concl} concludes the paper and presents a brief discussion about ongoing and future work.

%%%%%%%%%%%%%%%%%%%%%%%%%%%%%%%%%%%%%%%%%%%%%%%%%%%%%%%%%%%%%%%%%%%%%%%%%%%%%%%%
\section{MIXED-INTEGER QUADRATIC PROGRAMMING} \label{sec:MIQP}

We first introduce some of the basic concepts in mixed-integer programming based on branch-and-bound methods, such as convex relaxations and branching strategies.
%In direct optimal control, the continuous time OCP is first approximated by a tractable non-linear program. We formulate an equidistant grid over the control horizon $T$ consisting of the collection of time points $\{t_i\}^{N}_{i=0}$. We also consider a piecewise constant control parameterization $u(\tau) = u_i$ for $\tau \in [t_i,t_{i+1})$. 

\subsection{Convex Quadratic Program Relaxations}
\label{sec:qprel}

%The MIQP in~\eqref{MI-OCP} is a hard problem to solve in general. 
A standard approach to solve the MIQP~\eqref{MI-OCP} is to create convex relaxations of this problem (by either dropping some constraints or by re-formulating the problem and providing an approximation scheme) and then solve the relaxations in order to approach the solution to the original MIQP. A straightforward idea is to obtain convex QP relaxations by dropping the binary equality constraints~\eqref{MS:int} and instead enforcing the affine inequality constraints $0 \le F_iu_{i} \le 1$. Other convex relaxations for MIQPs have been studied in the literature such as moment or SDP relaxations that are often tighter than QP relaxations~\cite{luo2010semidefinite,Axehill2010}, but they can be relatively expensive to solve for larger problems. 

%In addition, several approximation schemes have been proposed for a great variety of MIQPs, such as randomization and graph-based approximations~\cite{lewin2002improved},~\cite{goemans2004approximation}, outer approximation schemes~\cite{fletcher1994solving} and simple heuristic based methods~\cite{takapoui2017simple}.
%\todo{add citation for graph-theory relaxations}

For the purpose of this paper, we will focus our attention on QP relaxations where we allow the binary variables to take on real values. The main reason for choosing this relaxation is that we utilize a tailored structure exploiting active-set solver, called \software{PRESAS}~\cite{PRESAS}, proposed recently for efficiently solving the convex QP relaxations. The latter solver has been shown to be competitive with state-of-the-art QP solvers for embedded MPC, and it benefits strongly from warm-starting, which can be exploited when solving the sequence of QPs within the branch-and-bound strategy. Note that the relaxations need to be convex, i.e., the weight matrices $Q_i$, $R_i$ and $P$ need to be positive (semi-) definite in~\eqref{MS:obj} such that each solution to a QP relaxation is globally optimal. 
%The performance of this solver when compared to other state of art solvers for discrete-time OCPs can be found here [cite].  We combine our active-set solver with the branch-and-bound method in order to solve the MIQP to optimality. The active-solver allows warm-starting which is particularly important when it is used embedded in the branch-and-bound method.

\subsection{Branch-and-Bound Algorithm}
\label{sec:BB_OCP}

\begin{figure} 
\centering
\includegraphics[trim={0in 1.3in 0.3in 0.5in},clip,scale=0.26]{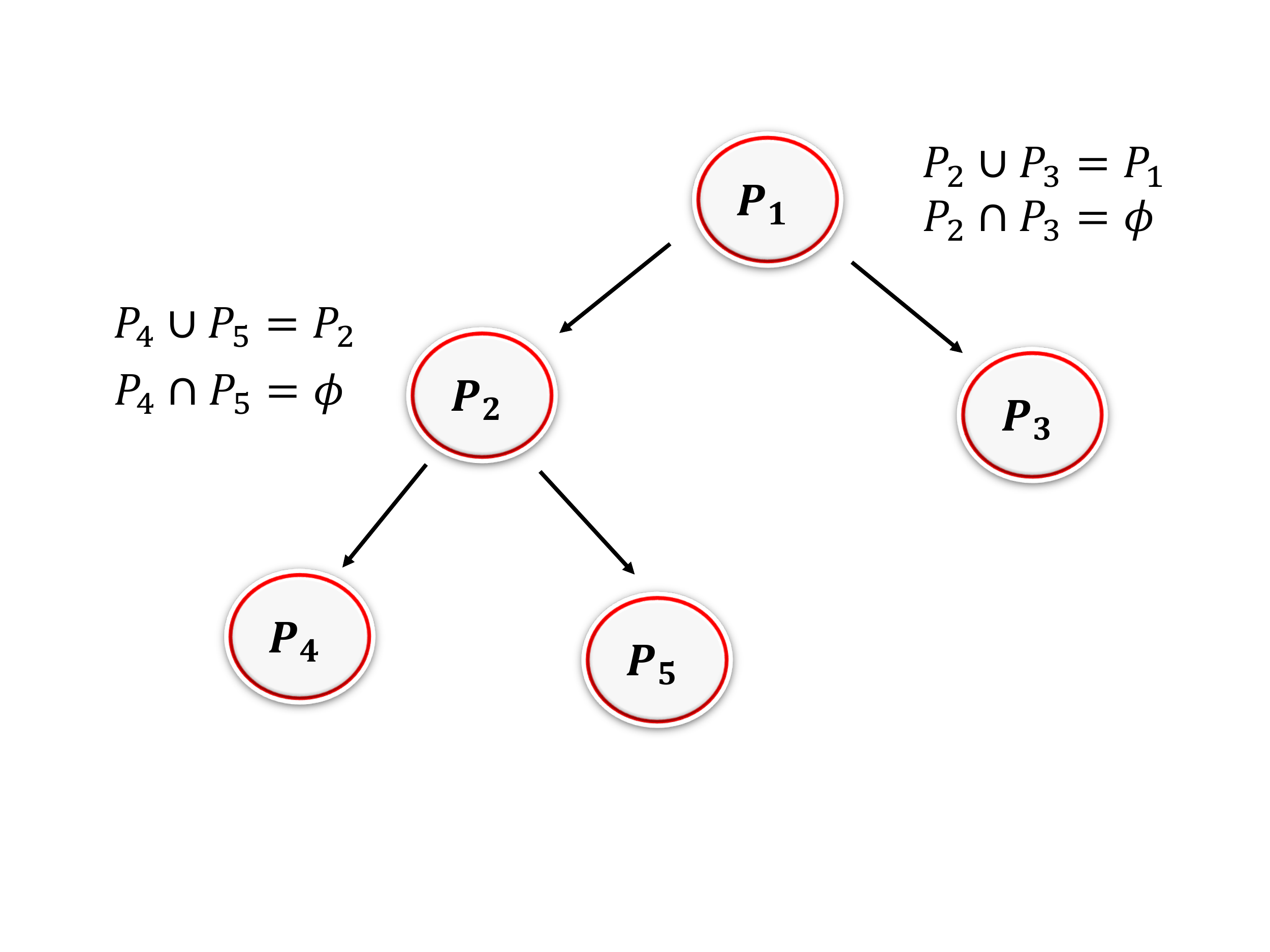}
\caption{Illustration of the branch-and-bound method as a binary search tree. A selected node can be either \emph{branched}, resulting in two partitions for each bound value in~\eqref{MS:int}, or \emph{pruned} based on feasibility or the current upper bound.}
\label{fig:BB_ilu}
\end{figure}

The main idea of the branch-and-bound~(B\&B) algorithm is to sequentially create partitions of the original problem and then attempt to solve those partitions. While solving each partition may still be challenging, it is fairly efficient to obtain local lower bounds on the optimal objective value, by solving relaxations of the mixed-integer program or by using duality. If we happen to obtain an integer-feasible solution while solving a relaxation, we can then use it to obtain a global upper bound for the solution to the original problem. This may help to avoid solving or branching certain partitions that were already created, i.e., these partitions or nodes can be \emph{pruned}. The general algorithmic idea of partitioning is better illustrated as a binary search tree, see Figure~\ref{fig:BB_ilu}.

A key step in this approach is how to create the partitions, i.e., which node to choose and which binary variable to select for branching. Since we solve a QP relaxation at every node of the tree, it is natural to branch on one of the binary variables with fractional values in the optimal solution of the QP relaxation. Therefore, if a variable, e.g., $u_{i,k} \in \{0, 1\}$ has a fractional value in a given QP relaxation, then we create two partitions where we respectively add the equality constraint $u_{i,k} = 0$ and $u_{i,k} = 1$.
%Therefore if a given variable, say $u_1$,} has a fractional value of $f$ in a given QP relaxation, then we create two sub-problems, 
Another key step is how to choose the order in which the created subproblems are solved. These two steps have been extensively explored in the literature and various heuristics are implemented in state-of-the-art tools~\cite{achterberg2005branching}. We provide next a brief description of strategies that we implemented in our B\&B solver.
%\todo{What did you mean by ``where we respectively add the constraints $u_1 \leq f -1$ to one problem and $u_1 \geq f +1$ to the other.''?}

\subsection{Tree Search: Node Selection Strategies}

A common implementation of the branch-and-bound method is based on a \emph{depth-first} node selection strategy, which can be readily implemented using a last-in-first-out~(LIFO) buffer. The next node to be solved is selected as one of the children of the current node and this process is repeated until a node is pruned, i.e., the node is either infeasible, optimal or dominated by the upper bound, which is followed by a backtracking procedure. Instead, a \emph{best-first} strategy selects the node with the lowest local lower bound so far. In what follows, we will employ a combination of the depth-first and best-first node selection approach. This idea is motivated by aiming to find an integer-feasible solution quickly at the start of the branch-and-bound procedure~(depth-first) to allow for early pruning, followed by a more greedy search for better feasible solutions~(best-first).

\subsection{Reliability Branching for Variable Selection}
%On this subsection we briefly describe the branching rules that will be used in our algorithm. The original idea was presented in [cite]. 
The idea of \emph{reliability branching} is to combine two powerful concepts for variable selection: strong branching and pseudo-costs~\cite{achterberg2005branching}. Strong branching relies on temporarily branching, both up (to higher integer) and down (to lower integer), for every binary variable that has a fractional value in the solution of a QP relaxation in a given node, before committing to the variable that provides the highest value for a particular score function. The increase in objective values $\Delta_{i,k}^+$, $\Delta_{i,k}^-$ are computed when branching the binary variable $u_{i,k}$, respectively, up and down. Given these quantities, a simple scoring function $\text{score}(\cdot,\cdot)$ is computed for each binary variable. For instance, based on the product~\cite{le2017abstract}:
\begin{equation}
S_{i,k} = \text{score}(\Delta_{i,k}^-,\Delta_{i,k}^+) = \text{max}(\Delta_{i,k}^+,\epsilon) \, \cdot \, \text{max}(\Delta_{i,k}^-,\epsilon), \label{eq:scorefun}
\end{equation}
given a small positive value $\epsilon > 0$.
This branching rule has been empirically shown to provide smaller search trees in practice~\cite{achterberg2005branching}. The downside is that this procedure is relatively expensive since several QP relaxations are solved in order to select one variable to branch on.

%\todo{In strong branching, for each variable $u^{k}_i$ that is fractional in the solution of the QP relaxation we solve an two QP's where we temporarily branch on $u^{k}_i$ and we compute the increments on the objective function value. Let $\Delta^{-}_{k,i}$ and $ \Delta^{+}_{k,i}$ be the change in objective value: 
%\begin{equation}
%\Delta^{-}_{k,i} = Z_{-} - Z_{0} \text{  ,  } \Delta^{+}_{k,i} = Z_{+} - Z_{0}
%\end{equation}
%
%where $Z_{0}$ is objective value of the current QP relaxation and $Z^{+},Z^{-}$ are the value of two temporary children. Then we can compute a score $S^{k}_{i} =\text{score}(\Delta^{-}_{k,i}, \Delta^{+}_{k,i}) $, for some scoring function (for example the product). We repeat this procedure for every integer decision variable that takes a fractional value in the optimal solution at the current node. The chosen branched variable is the one that has he highest score.}

The idea of pseudo-costs aims at approximating the increase of the objective function to decide which variable to branch on, without having to solve additional QP relaxations. This can be done by keeping statistic information for each binary variable, i.e., the \emph{pseudo-costs} that represent the average increase in the objective value per unit change in that particular binary variable when branching. Every time that a given variable is chosen to be branched on, and the resulting relaxation is feasible, then we update each corresponding pseudo-cost with the observed increase in the objective, divided by the distance of the real to the binary value, in the form of a cumulative average. Therefore, each variable has two pseudo-costs, $\phi_{i,k}^{-}$ when the variable was branched ``down'' and $\phi_{i,k}^{+}$ when it was branched ``up''. Given the solution to a QP relaxation, one can then use the pseudo-costs to select the binary variable with the highest score value to be branched on next:
\begin{equation}
S_{i,k} = \text{score}(\bar{u}_{i,k} \,\phi_{i,k}^-, \;(1-\bar{u}_{i,k}) \,\phi_{i,k}^+), \label{eq:score_pc}
\end{equation}
given a fractional value $\bar{u}_{i,k}$ in the QP relaxation.
%In order to decide if on a given subproblem we should branch on a given variable we utilize a score function, which is typically the product of both pseudo-costs. 

This way, we select variables based on their past behavior throughout the branch-and-bound tree. However, at the beginning of the algorithm, the pseudo-costs are not yet initialized, which is when branching decisions typically impact the tree size the most. \emph{Reliability branching} uses strong branching to initialize the pseudo-costs until a certain condition of reliability is satisfied, e.g., one switches to using pseudo-costs only once that particular variable has been branched on a specified number $\eta_{rel}$ of times~\cite{achterberg2005branching}. 
%This number is called the reliability threshold $\eta_{rel}$. 
The resulting branching rule is summarized in Algorithm~\ref{RelBranch}. Note that reliability branching coincides with pseudo-cost branching if $\eta_{rel}=0$, with strong branching if $\eta_{rel}=\infty$, but typically a value $1 \le \eta_{rel} \le 4$ is chosen.
%\todo{Suppose that at a previously explored node, say Q', we branched on variable $u^{k}_i$. Again we let $Z_{-}$ and $Z_{+}$ be the respective children objective values. Let $f^{k}_i$ be the fractional value of variable $u^{k}_i$ in the optimal solution at Q'. Then we can compute the relative increase $\psi^{-}_{k,i},\psi^{+}_{k,i}$:
%\begin{equation}
%\psi^{-}_{k,i} = \frac{Z_{-} - Z_{0}}{f^{k}_i} \text{  ,  } \psi^{+}_{k,i} = \frac{Z_{+} - Z_{0}}{1 - f^{k}_i}
%\end{equation}
%
%Now let $\sigma^{-}_{k,i},\sigma^{+}_{k,i}$ be the sum of $\psi^{-}_{k,i},\psi^{+}_{k,i}$, respectively over all nodes where variable $u^{k}_i$ was selected to be branched on. Let $\eta^{-}_{k,i},\eta^{+}_{k,i}$ be the number of problems that were feasible upon branching on $u^{k}_i$ downwards and upwards, respectively. Then we can define the pseudo-costs as:
%\begin{equation}
%\phi^{-}_{k,i} = \frac{\sigma^{-}_{k,i}}{\eta^{-}_{k,i}} \text{  ,  } \phi^{+}_{k,i} = %\frac{\sigma^{+}_{k,i}}{\eta^{+}_{k,i}}
%\end{equation}
%
%Now we can compute the score of variable $u^{k}_i$ as $S^{k}_{i} =\text{score}(f^{k}_i\phi^{-}_{k,i}, (1-f^{k}_i)\phi^{+}_{k,i}) $.
%Reliability branching combine both ideas: If a given decision variable was not branched on a certain number of times, then we perform strong branching on it, otherwise we call that variable reliability, and we can fore-go the need to solve the relaxations and utilize the pseudo-costs to compute that variable score. The whole reliability branching scheme is summarized below:}

This rule can be further augmented by implementing a look ahead limit in the number of candidates, as well as a limit in the number of iterations for each QP relaxation in the strong branching step. 
%We refer to [ref] for a more detailed description of those branching rules. 
Note that many other branching rules exist such as, e.g., ``most infeasible'' branching which selects the binary variable with fractional part that is closest to $0.5$. Even though the latter rule is used quite often, e.g., in~\cite{Bemporad2018}, it generally does not perform very well in practice~\cite{achterberg2005branching}.
Extensive empirical experiments with different branching strategies are beyond the scope of this paper.
%will be present in th journal version of this paper.

    \begin{algorithm}[tpb]
    	\caption{Reliability Branching Strategy}
    	\label{RelBranch}
    	\begin{algorithmic}[1]
    		\Require $\eta_{\text{rel}}$, set C of candidate variables for branching. 
%    		Reliability threshold $\eta_{\text{rel}}$.
    		\For {candidate variables $u_{i,k}$ in C}
    		
    		\If {$\#\texttt{branch}(u_{i,k}) \le \eta_{\text{rel}} $}
    		\State Strong branching on $u_{i,k}$ to compute score $S_{i,k}$. 
    		\State Update pseudo-costs $\phi^{-}_{i,k}$ and $\phi^{+}_{i,k}$.
    		\Else
    		\State $S_{i,k} = \text{score}(\bar{u}_{i,k} \,\phi_{i,k}^-, \;(1-\bar{u}_{i,k}) \,\phi_{i,k}^+)$.
    		\EndIf
    		\EndFor
    		\Ensure Select variable with highest score $S^{*} = \underset{i,k}{\max}\,{S_{i,k}}$.
    	\end{algorithmic}
    \end{algorithm}
    
%%%%%%%%%%%%%%%%%%%%%%%%%%%%%%%%%%%%%%%%%%%%%%%%%%%%%%%%%%%%%%%%%%%%%%%%%%%%%%%%
\section{PRESOLVE TECHNIQUES FOR MIXED-INTEGER OPTIMAL CONTROL} \label{sec:PRESOLVE}

As mentioned earlier, presolve techniques are often crucial in making convex relaxations tighter such that typically fewer nodes need to be explored, sometimes to such an extent that seemingly intractable problems become tractable. Next, we briefly describe some of these concepts with a focus on domain propagation for bound strengthening and its implementation for mixed-integer optimal control.

\subsection{Domain Propagation for Condensed QP Subproblem}

%Next to branching and node selection strategies, a fundamental way of improving the performance of the branch-and-bound algorithm is by providing stronger relaxations, resulting in higher local lower bounds for the partitions. 
Several strengthening techniques are implemented as part of ``presolve'' routines in commercial solvers~\cite{achterberg2016presolve}. One particular technique that is suitable to mixed-integer optimal control is based on \emph{domain propagation}, in which the goal is to strengthen bound values based on the inequality constraints~\eqref{MS:path}-\eqref{MS:term} in the problem. However, the results of such a strategy are rather weak when directly applied to the block-sparse QP in~\eqref{MI-OCP}, because the stage-wise coupling of the state variables~\eqref{MS:Kstep} needs to be taken into account. Therefore, we use instead the equivalent dense QP formulation in which the state variables are numerically eliminated, such that stronger bounds can be obtained for the control variables. Hence, we can use the block-structured sparsity to efficiently solve the QP relaxations, while we use the equivalent but dense format to effectively perform domain propagation.

%We will leverage the fact that we can eliminate the states variables to obtain a dense QP problem (a process we refer to as Condensing). The resulting problem has a dense Hessian and constraint matrices. That allow us to derive stronger bounds for the control variables, which are still valid for the original formulation. Hence we can use the OCP structure to not only obtain tighter bounds (via the condensed formulation) but we can keep the original formulation when we solve the QP relaxations (by using a sparsity exploiting-solver).

Let us concatenate all state variables in a vector $X$ and all control variables in the vector $U$, such that Eqs.~\eqref{MS:initial}-\eqref{MS:Kstep} can be written more compactly as
\begin{equation}
\bar{A}X = \bar{B}U + b + {E_0}\hat{x}_0,
\end{equation}
where we define the block-sparse matrices 
 \begin{subequations}
 \begin{alignat}{2}
 \bar{A} &= \begin{bmatrix}
 I & & & \\
 -A_{1} & I & & \\
%  & -A_{2} & I & \\
  &  \ddots & \ddots \\
  & & -A_{N-1} & I
 \end{bmatrix},  \\
 \bar{B} &= \text{blkdiag}(B_0,\ldots,B_{N-1}) \text{   ,   } E_0 = [A^\top_0, 0, \ldots, 0]^\top. 
 \end{alignat}
 \end{subequations}
The matrix $\bar{A}$ is invertible such that we can write:
 \begin{equation}
X = \bar{A}^{-1}\bar{B}U + \bar{A}^{-1}(b + {E_0}\hat{x}_0).
 \end{equation}
Now, we can substitute the latter expression for the state vector in OCP~\eqref{MI-OCP} to obtain the condensed form
\begin{subequations} \label{Cond-OCP}
	\begin{alignat}{5}
	\underset{U}{\text{min}} \quad &  \frac{1}{2} U^{\top}\Hc U + \hc^{\top} U&& \label{cond:obj}\\
	\text{s.t.} \quad\; & \bar{l}^{c} \leq \Dc U \leq \bar{u}^{c} \label{cond:path} \\
	& F_iu_{i} \in \{0, 1\}, \quad && i \in \{0,\ldots,N-1\}, \label{cond:int} 
%	&  u^{k}_{i} \in \mathbb{Z} \quad &&k \in I_{k}, i \in [N-1] \label{cond:setU} 
	\end{alignat}
\end{subequations}
where the condensed matrices and vectors read as
 \begin{subequations}
 \begin{alignat}{4}
 \Hc &= (\bar{A}^{-1}\bar{B})^{\top}Q\bar{A}^{-1}\bar{B} + R, \; \Dc = C\bar{A}^{-1}\bar{B} + D, \\
 \hc &= (\bar{A}^{-1}\bar{b})^{\top}Q\bar{A}^{-1}\bar{B}, \\
  \bar{l}^{c} &= l^{c} - C\bar{A}^{-1}\bar{b}, \quad \bar{u}^{c} = u^{c} - C\bar{A}^{-1}\bar{b}, 
 \end{alignat}
 \end{subequations}
where $\bar{b} := b + E_0\hat{x}_0$ is defined and given $Q = \text{blkdiag}(Q_{1},\ldots, Q_{N-1}, P)$ , $R = \text{blkdiag}(R_{0},\ldots, R_{N-1})$, and $l^{c} = [l^{c^\top}_{1}, \ldots, l^{c^\top}_{N}]^{\top}$ and $u^{c} = [u^{c^\top}_{1}, \ldots, u^{c^\top}_{N}]^{\top}$. 

Given the condensed problem formulation, which can be computed offline and which is parametric in the current state value $\hat{x}_0$, we can then apply the following bound strengthening procedure, which is explained next for a single affine constraint $l_b \le \sum_{i} d_i u_i \leq u_b$ in~\eqref{cond:path}. This constraint can be used to try and tighten bound values for all control variables $u_i$ for which $d_i \neq 0$, where $u_i$ denotes a single control variable in the vector $U$. 
%As an illustration we show it below for a single constraint of the type $\{d^{\top}U \leq u\}$ and for a variable $u_1$. 
Let $\bar{u}_{i}, \underline{u}_{i}$ be the current upper/lower bounds for $u_{i}$ such that
\begin{equation}
d_i u_i \leq u_b - \sum_{j \neq i} d_{j}u_{j} \leq \underbrace{u_b - \sum_{j \neq i, d_{j} > 0} d_{j}\underline{u}_{j} - \sum_{j \neq i, d_{j} < 0} d_{j}\bar{u}_{j}}_{=: \bar{u}_{b,i}},
\end{equation}
in which we divide by $d_i$ in order to obtain
\begin{equation}
u_i \leq  \frac{\bar{u}_{b,i}}{d_i} , \text{ if } d_i > 0 \quad \text{or} \quad 
u_i \geq \frac{\bar{u}_{b,i}}{d_i} , \text{ if } d_i < 0. \label{DPEQ}
\end{equation}
This results, respectively, in the updated bound values
\begin{equation}
\bar{u}_{i} = \min{(\bar{u}_{i},\frac{\bar{u}_{b,i}}{d_i})} , \quad \text{or} \quad 
\underline{u}_{i} = \max{(\underline{u}_{i},\frac{\bar{u}_{b,i}}{d_i})},  \label{eq:bound_up}
\end{equation}
or, in case $u_i$ is an integer or binary variable,
\begin{equation}
\bar{u}_{i} = \min{(\bar{u}_{i},\left\lfloor\frac{\bar{u}_{b,i}}{d_i}\right\rfloor)} , \quad \text{or} \quad 
\underline{u}_{i} = \max{(\underline{u}_{i},\left\lceil\frac{\bar{u}_{b,i}}{d_i}\right\rceil)}. \label{eq:bound_up2}
\end{equation}

where $\lfloor \cdot \rfloor$ and $\lceil \cdot \rceil$ are the floor and ceiling operations, respectively. Thus, this can result in strengthening of bound values for both continuous and integer/binary control variables. The procedure can be executed for each control variable and each inequality constraint in an iterative manner, see Algorithm~\ref{DomainPropagation}, since bound strengthening for one variable can lead to strengthening for other variables~\cite{achterberg2016presolve}. The process is typically stopped when the bound values do not sufficiently change or a certain limit on the computation time is met.

Domain propagation can lead to considerable reductions in the amount of explored nodes, e.g., because variables are fixed, when $\bar{u}_{i} = \underline{u}_{i}$, or because of infeasibility detection, when $\bar{u}_{i} < \underline{u}_{i}$, without the need to solve any QP relaxations. In addition, the updated bound values for all control variables can be used to strengthen QP relaxations in the future.
Lastly, we can use domain propagation in order to improve and generalize Hessian-based fixing strategies, such as the one proposed in~\cite{Axehill2004}. Hessian-based fixing typically can only be applied to unconstrained problems, since it fixes the variables solely based on the objective. Here, we propose to use domain propagation to compute the feasibility impact of certain variable fixings. More specifically, a particular variable can be fixed based on optimality, if and only if this fixing does not induce feasibility-based fixings.
%More specifically, a particular variable can be fixed based on optimality, if and only if this fixing does not lead to other feasibility-based fixings such that no integer-feasible solutions are excluded.

%If the variable in question is constrained to take integer values we can take the floor and ceil, respectively to improve even more the bound obtained. 
%The above computation can be repeated for every constraint and variable. The whole procedure is summarized in Algorithm~\ref{DomainPropagation}. 
%Note that while during this procedure if $\bar{l}_{i} > \bar{u}_{i}$, then we can safely declare the relaxation infeasible, without even solving it. So this procedure can also result in a certificate of infeasibility for the QP relaxation.

\begin{algorithm}[tpb]
    	\caption{Domain Propagation for Bound Strengthening}
    	\label{DomainPropagation}
    	\begin{algorithmic}[1]
    		\Require Inequality constraints~\eqref{cond:path}, variable bounds $\bar{u}_{i}, \underline{u}_{i}$.
%    		Condensed constraint matrix $\bar{D}$, LHS and RHS vectors $\bar{l}^{c}, \bar{u}^{c}$. Current variable bounds  $\bar{l}, \bar{u}$.
			\While {stopping criterion == False}
			\For { every row of $\Dc$}
    		\For { every $u_{i} \in U, \;d_i \neq 0$}
    		\State Obtain bound values $\bar{u}_{b,i}, \bar{l}_{b,i}$ using Eq.~\eqref{DPEQ}.    	
    		\State Update variable bounds using~\eqref{eq:bound_up} or~\eqref{eq:bound_up2}.
%    		\State Update the bounds $\bar{u}_{i} \leftarrow \min(\bar{u}_{i}, \min_{j}(ub_{j}) ) $
%    		\State Update the bounds $\bar{l}_{i} \leftarrow \max(\bar{l}_{i}, \max_{j}(lb_{j}) ) $
    		\EndFor
    		\EndFor
    		\EndWhile
			\Ensure Updated bounds $\bar{u}_{i}, \underline{u}_{i}$ for all control variables.
    	\end{algorithmic}
    \end{algorithm}

%\subsection{Hessian based Variable Fixings}

%\todo{brief description here as well or refer this to future work, e.g., the journal extension?}
%\todo{cite~\cite{Axehill2004} for this work, which we extended to the constrained case by using domain propagation to check whether fixings can be made or not}

%\subsection{Dual Active-Set Iterations}
%
%\todo{brief description here as well?}
%\todo{cite~\cite{Fletcher1998} for this}
%\todo{Maybe leave out postsolve strategies for now?}
    
\subsection{Probing Strategies and Cutting Planes}

Probing~\cite{savelsbergh1994preprocessing} is a classical technique that can be incorporated in any branch-and-bound method to derive stronger inequalities or better bounds. It consists of tentatively trying to fix some variables and to derive potential logical implications on other variables. We do not further describe probing strategies in detail, but we refer to~\cite{achterberg2016presolve} for an overview. 
%The typical trade-off in designing probing strategies is that, while it can lead to additional fixings of the integer variables, a considerable amount of additional computations are often required. 
The computational cost and performance of probing can be greatly improved by relying on some of the other techniques that were discussed earlier. For example, the pseudo-costs can be used in order to choose the bound value for each binary variable that is likely to result in a low objective value. In turn, the QP relaxations that are solved in the probing procedure can be used to update the pseudo-cost values. In addition, domain propagation and other variable fixing strategies can be used to reduce the amount of QP relaxations that need to be solved.
	
%Therefore, this technique can not be extensively applied throughout the algorithm. However, this heuristic strategy can be potentially improved by using information about the variables such as their bounds and the pseudo-costs.

%Lastly, we will not explore cut generation strategies, although any cut generated for the condensed problem can be fully transferred to the original OCP formulation. The reason for that is without a careful and tailored cut generation procedure, the inequalities would couple variables, potentially across stages, which is undesirable as we rely on a block-sparsity exploiting QP solver.
Other presolving techniques such as cut generation can be applied using the condensed problem, and can be fully transferred to the original OCP formulation. In the present paper, we refrain from using cut generation techniques as they produce inequalities that potentially couple variables across stages. Such coupling between stages is not desirable as we rely on a block-sparsity exploiting QP solver.

%%%%%%%%%%%%%%%%%%%%%%%%%%%%%%%%%%%%%%%%%%%%%%%%%%%%%%%%%%%%%%%%%%%%%%%%%%%%%%%%%
%\section{PROPOSED MIQP ALGORITHM FOR OPTIMAL CONTROL}
\subsection{Resulting MIQP Algorithm for Optimal Control}

Algorithm~\ref{MIQP_ALG} describes the most important steps in our proposed B\&B method for solving the MIQP in~\eqref{MI-OCP}. It solves a block-structured QP relaxation using \software{PRESAS}~\cite{PRESAS} at every node and utilizes reliability branching (Algorithm~\ref{RelBranch}) to decide the branching variables. As discussed earlier, the node selection strategy is based on a depth-first search followed by a best-first search as soon as an integer-feasible solution has been found. Note that the upper bound value UB provided to Alg.~\ref{MIQP_ALG} can be based on an integer-feasible solution guess or it can be set to $+\infty$.
%The Node Selection Strategy presented in the Algorithm is the Depth-First Search (LIFO: last in, first out), as the new generated subproblems are appended at the end of list and are explored next. It is typical to start of a LIFO Node Selection Strategy until a feasible solution is find and then switch to a Best Bound approach (BESTBOUND), where we choose the node which contains the most promising lower-bounds (smallest). The idea is that once feasibility is ascertained we can then focus on seeking optimality. 
Because of space limitations, the present paper will not further discuss all parameter choices in the algorithm such as, e.g., the reliability branching parameters, presolving frequency, memory usage, etc. 
%The study of the impact of all different design parameters in the performance of the our MIQP solver and it's impact on the application for MPC is part of ongoing work.

  \begin{algorithm}[tpb]
    	\caption{B\&B Method for the MIQP-OCP in~\eqref{MI-OCP}}
    	\label{MIQP_ALG}
    	\begin{algorithmic}[1]
    		\Require Upper bound UB, tolerance $\epsilon$.
    		%\Statex \texttt{Problem linearization}
    		\State LB=$-\infty$ and initialize $L=\{P_0\}$ with root node.
    		\State Select current node $P_{c} \gets P_0$.
    		\While{$UB-LB > \epsilon$}
    		\State{Apply domain propagation to $P_c$ using Alg.~\ref{DomainPropagation}.}
    		\State{Solve resulting QP relaxation with \software{PRESAS}.}
    		
    		\If {QP is feasible and $J(\bar{X},\bar{U}) \leq $ UB}

    		\If {QP solution is not integer-feasible}
    		\State $LB \gets \min_{P \in L}{J(P)}$.
    		\State  {Select branching variable $v$ using Alg.~\ref{RelBranch}.} 
    		\State  {Create subproblems $P_u$~``up'' and $P_l$~``down''.}
%    		\State {v $\leftarrow$ Reliability Branching (Algorithm 1).}
    		\State Append $\{P_l, P_u\}$ to $L$ if $(1-\bar{v})\phi^{+}_{v} < \bar{v} \phi^{-}_{v}$
    		\Statex \hspace{4.5em}{or append $\{P_u, P_l\}$ to $L$, otherwise.}
    		
%    		\State  Append the formed subproblems $P_l$ and $P_u$ to the list $L$.
%    		\State {Remove $P_c$ from $L$.}
    		\Else \State $UB \gets J(\bar{X},\bar{U})$ and $(X^{*},U^{*}) \gets (\bar{X},\bar{U})$.
%    		\State Store the solution vector $(x^{*},u^{*})$ found on node $P_c$ as the incumbent solution.
    		\EndIf
%    		\Else \State {Remove $P_c$ from $L$.}.   		
    		\EndIf
    		\State {Remove current node $P_c$ from to-do list in $L$.}
    		\State {Select next node based on depth-first (last node}
    		\Statex \hspace{1.5em}{in list $L$) or based on best lower bound.}
%    		\State LIFO if no feasible has been found, BESTBOUND otherwise.
    		\EndWhile
    		\Ensure MIQP solution vector $(X^{*},U^{*})$.
    	\end{algorithmic}
    \end{algorithm}

%%%%%%%%%%%%%%%%%%%%%%%%%%%%%%%%%%%%%%%%%%%%%%%%%%%%%%%%%%%%%%%%%%%%%%%%%%%%%%%%
\section{MIXED-INTEGER MPC ALGORITHM}
\label{sec:three}
In embedded applications of mixed-integer MPC, one needs to solve an MIQP~\eqref{MI-OCP} at each sampling instant under strict timing constraints. We can leverage the fact that we solve a sequence of similar problems (parametrized by the initial condition $\hat{x}_{0}$), in order to warm-start the B\&B-algorithm. We refer to our proposed warm-starting procedure as \emph{tree propagation}, because the main goal is to ``propagate'' the B\&B tree forward by one time step. We describe this process in detail below. Then, we present the resulting mixed-integer MPC algorithm. 

\subsection{Warm Starting based on Tree Propagation}

\begin{figure} 
\centering
\includegraphics[trim={0.1in 1.0in 0.2in 0.7in},clip,scale=0.26]{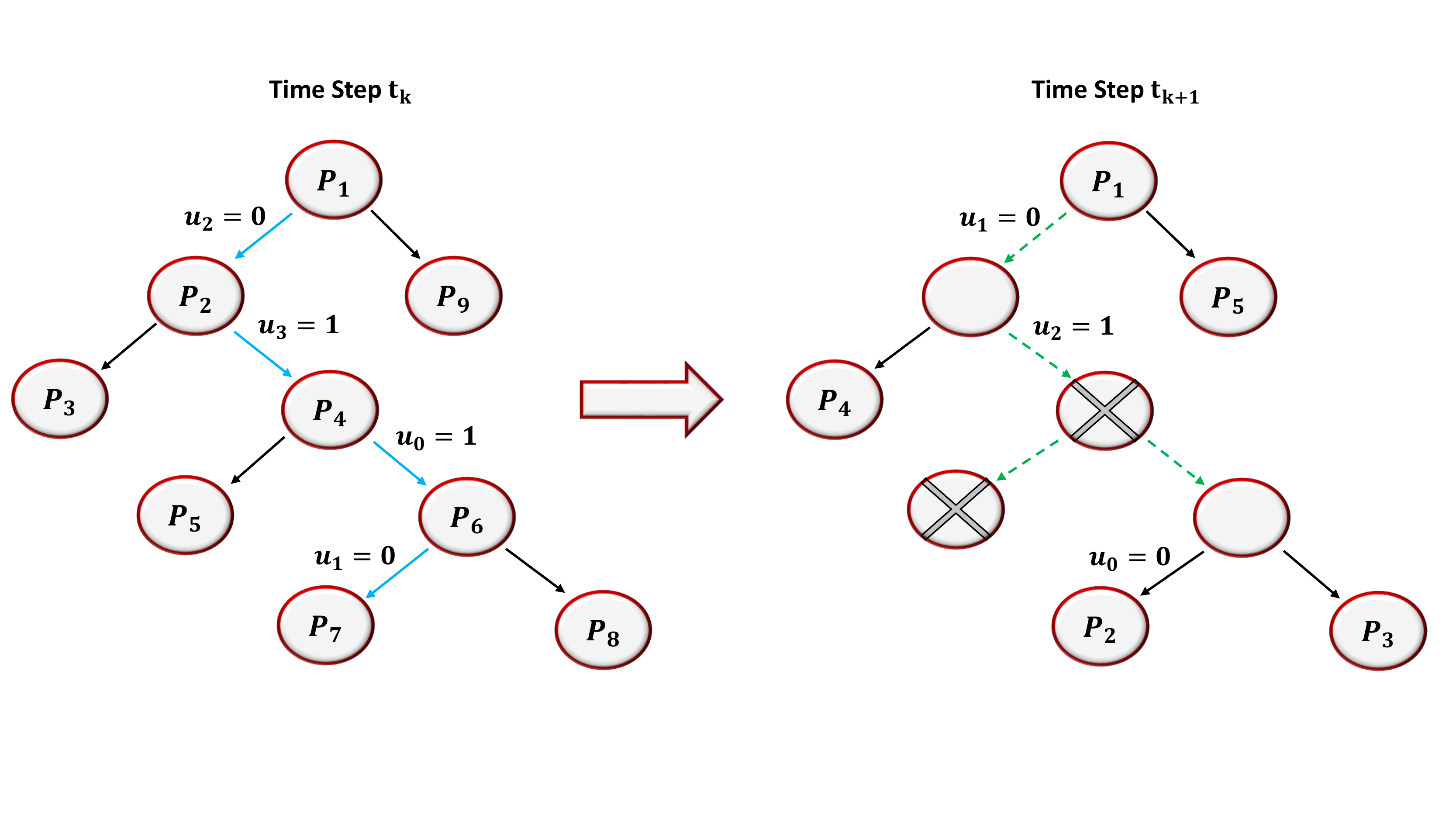}
\caption{Illustration of the tree propagation technique from one time point to the next in the MI-MPC algorithm: index~$i$ denotes the order in which each node $P_i$ is solved.}
\label{fig:MPCTree}
\end{figure}

The warm-starting procedure aims to use knowledge of one MIQP, i.e., the search tree after solving the problem, in order to improve the B\&B search for the next MIQP. Our idea is to store the path from the root to the leaf node where the optimal solution to the MIQP was found, as well as the branching order of the variables. We can then perform a shifting of this path in order to obtain a ``warm-started tree'' to start our search to solve the MIQP at the next time step. We illustrate this procedure in Figure~\ref{fig:MPCTree}, where the optimal path at the current time step is denoted by the sequence of nodes $P_1 \rightarrow P_2 \rightarrow P_4 \rightarrow P_6 \rightarrow P_7$. Let us consider a corresponding sequence of variables $u_2 \rightarrow u_3 \rightarrow u_0 \rightarrow u_1$ that we branched on in order to create such optimal path. After shifting by one time step, all branched variables in the first control interval can be ignored, e.g., resulting in a shifted and shorter path of variables $u_1 \rightarrow u_2 \rightarrow u_0$.    
%On the figure, at stage the optimal path (i.e.: the branch where we found the optimal solution) is $P_1 \rightarrow P_2 \rightarrow P_5 \rightarrow P_6 \rightarrow P_7$. 
%And the variables we branched on were $u_2 \rightarrow u_3 \rightarrow u_0 \rightarrow u_1$. Now we shift the indices of the branched variables by setting $i \leftarrow i-1$. But, node that is not possible for $u_0$ as it is the initial stage, hence we prune this link from the shifted path. Hence the shifted path for the branched variables becomes $u_1 \rightarrow u_2 \rightarrow u_0$. 

At the subsequent time step, after obtaining the new state estimate, we execute all presolving techniques and we solve the QP relaxation corresponding to the root node. After removing from the warm-started tree the nodes that correspond to branched variables which are already integer feasible in the relaxed solution at the root node, we proceed by solving all the leaf nodes on the warm-started path. As we solve both children of a node on this path, we do not have to solve the parent node itself and therefore reduce computations by solving less QP relaxations.
%as always ($P_1$), but now instead of solving the direct children of the root node, we proceed to directly to solve the leaf of warm-started branch. We then proceed in a depth-first search manner. But note that as we solve both children of the node we do not need to solve the parent. 
Hence, we go over the tree in the order depicted by the index of each node in Fig.~\ref{fig:MPCTree}. After the warm-started branch has been explored, we resume normal procedure of the B\&B method. Algorithm~\ref{alg:Tree_prop} summarizes the proposed tree propagation technique.

\begin{algorithm}[h]
	\caption{Tree Propagation for Warm-Started B\&B}
	\label{alg:Tree_prop}
	\begin{algorithmic}[1]
    		\Require Optimal path $P$ from root to leaf node.
    		\State Shift index of branched variables by $1$ stage along path.
    		\State Solve root node of shifted path $P$, including presolve.
    		\For {(branched variables on stage $-1$ after shifting) \newline
    			\hspace{1em}\OR~(variables are integer feasible in root node)\newline
    			\hspace{1em}\OR~(variables without pseudo-costs)}
    		\State {remove associated node from the path $P$.}  		
    		\EndFor	
    		\State Shift the QP relaxation solution on every node of the path and store it as a warm start for the QP solver.
    		\State Re-order sequence of branched variables by scoring based on warm-started pseudo-cost information.
    		\State Initialize the B\&B tree along the shifted path $P$, creating nodes along the path and their respective children.
    		\State Create the warm-started list $L$, excluding parent nodes.
    		\Ensure Warm-started tree for next MIQP, given by list $L$.		
	\end{algorithmic}
\end{algorithm}

The sequential nature of the problem also allows to shift and re-use the pseudo-cost information from one MPC time step to the next. This idea has the potential of producing smaller search trees as the MPC progresses, without the need to perform strong branching at every MPC step. The propagation of pseudo-costs can be coupled with an update of the reliability parameters to improve the overall performance. For example, the reliability number should be reduced for each variable from one time step to the next, in order to force strong branching for variables that have not been branched on in a sufficiently long time. In addition, nodes can be removed from the warm-started path in case they correspond to branched variables for which there is no pseudo-cost information or it is not sufficiently reliable, in an attempt to avoid bad branching decisions. Finally, these warm-started pseudo-costs can also be used to re-order the warm-started tree, in order to result in smaller search tree sizes.
%fine tuning of the Reliability parameters in order to make updates to pseudo-costs that are not being updated regularly and so on. 
%We highlight the usefulness of propagating pseudo-costs in our experiments.

%\todo{Discuss why we remove certain branched variables in the warm-started path of Algorithm~\ref{} + the use of pseudo-cost to ensure a good ordering of the branched variables. Essentially, the intuition is that this avoids situations where warm-starting is actually worse than cold-starting because you could branch on variables in the wrong order, or you could branch on variables that are already integer feasible because of presolving etc.}

The proposed tree propagation technique, with the additional re-use of pseudo-cost information, has been summarized in Algorithm~\ref{alg:Tree_prop}. This procedure can improve the overall performance of the B\&B method in multiple ways. First of all, the optimal path and pseudo-cost information is re-used to make better branching decisions for the mixed-integer program at the next time step, because the search trees are often similar for two subsequent problems. Also, the computational cost can be reduced by solving less QP relaxations to explore the warm-started tree. In addition, the shifted optimal path can be used in an attempt to efficiently obtain an integer-feasible solution, and therefore an important upper bound in the B\&B algorithm, for the MPC problem at the next time step. 
%We save some computation, as we do not need to solve the nodes of the warm-started path. Moreover, the B\&B tree is initialized in a sequence that is similar to the optimal path of the previous stage, so potentially it can be very similar for the new problem, thus resulting in fewer nodes to be explored. 
Lastly, one can store the relaxed QP solutions on the optimal path, shift them by one time step and use them to warm-start the QP solver for nodes on the shifted optimal path. 
%That aspect is very important as the QP solver may start from a better active set than a simple cold-start.

\subsection{MI-MPC Algorithm Implementation}

Algorithm~\ref{alg:MI_MPC} summarizes the proposed MI-MPC algorithm. It solves a sequence of MIQPs where the branch-and-bound tree is warm-started at every time step, as well as the pseudo-cost and QP condensing information. As mentioned earlier, the B\&B strategy and the additional presolve, warm-start and heuristic branching techniques have been implemented in \matlab, based on a C~code implementation of the \software{PRESAS} algorithm~\cite{PRESAS} to solve each QP relaxation. 
%The Condensed matrices are updated as the new initial state is obtained via Eq. 10.  We also shift the Pseudo-Costs information accordingly to be  re-used in the next stage's MIQP. The goal of this approach is two-fold: By re-using past information, such as branching, the hope is that they are still relevant in order to produce smaller search trees. On the other hand, be using the shifted optimal solution stored in a given node we are essentially warm-starting the active-set solver, that leads to essential speed-ups in the solution of the QP relaxations. Lastly, by warm-starting the tree the goal is to have a sequence of trees of similar size, hence decreasing the chance that "outliers" (that is a very large tree in a given time stage) to appear. 
In Section~\ref{sec:caseStudies}, we illustrate the computational performance of the presented MI-MPC algorithm, including these presolving and warm-starting techniques, for two numerical case studies of mixed-integer MPC. A self-contained C~code implementation is part of ongoing work, in order to illustrate the computational efficiency of the proposed algorithmic techniques.

\begin{algorithm}[h]
	\caption{Warm-Started B\&B Algorithm for MI-MPC}
	\label{alg:MI_MPC}
	\begin{algorithmic}[1]
    		\Require Current state $\hat{x}_0$, list of nodes $L$ and pseudo-costs.
    		\Statex \texttt{Solve MIQP}
    		\State Update condensing information, given current state $\hat{x}_0$.
    		\State Formulate MIQP~\eqref{MI-OCP} and solve it using Algorithm~\ref{MIQP_ALG}.	
    		\State Apply new control input $u_0^\star$ to the system. \vspace{1mm}
    		\Statex \texttt{Propagation Step}
    		\State Warm-start and shift pseudo-cost information.
    		\State Perform tree propagation to warm-start node list $L$ for the next MI-MPC time step (see Algorithm~\ref{alg:Tree_prop}).
%    		\Ensure MPC closed-loop controller $u_0^+$ and new system state  $x^{+}_{0}$. 
	\end{algorithmic}
\end{algorithm}

\subsection{Real-Time Embedded Applications of MI-MPC}

Note that, in practice, the proposed warm-starting strategies often allow one to obtain an integer-feasible solution in a computationally efficient manner. However, even if the tree propagation immediately provides the globally optimal solution to the MIQP~\eqref{MI-OCP}, a branch-and-bound algorithm still needs to perform relatively many iterations to prove optimality by pruning remaining nodes in the search tree. This motivates the use of a maximum number of B\&B iterations in order to meet strict timing requirements of the embedded control application. Even if the algorithm does not terminate within this specified number of iterations, a feasible or even optimal solution may be available. This and other heuristic strategies for real-time implementation of MI-MPC are straightforward but outside the scope of this paper.
%In what follows, we therefore solve each MIQP exactly.

\section{CASE STUDIES: MIXED-INTEGER MPC}
\label{sec:caseStudies}

We report two numerical case studies to illustrate the computational performance of our MIQP-based MPC algorithm: a hybrid MPC test example and a satellite orbit re-centering application with a no-go zone in the orbital path. Our branch-and-bound algorithm has been implemented in \matlab~in conjunction with the \software{PRESAS} active-set solver in C. To evaluate the performance, we compare our algorithm with the state-of-the-art \software{GUROBI}~\cite{gurobi} and \software{MOSEK}~\cite{mosek} solvers for mixed-integer programming. 
%In our comparisons, we replaced the entire branch-and-bound algorithm with \software{GUROBI} solver. 
It is important to emphasize that all advanced presolve and heuristic options have been activated for both software tools, resulting in fair computational comparisons.

\subsection{Hybrid MPC: Benchmark Example}
\begin{figure}[tpb] 
	\centerline{\hbox{
			\includegraphics[width=0.55\textwidth]{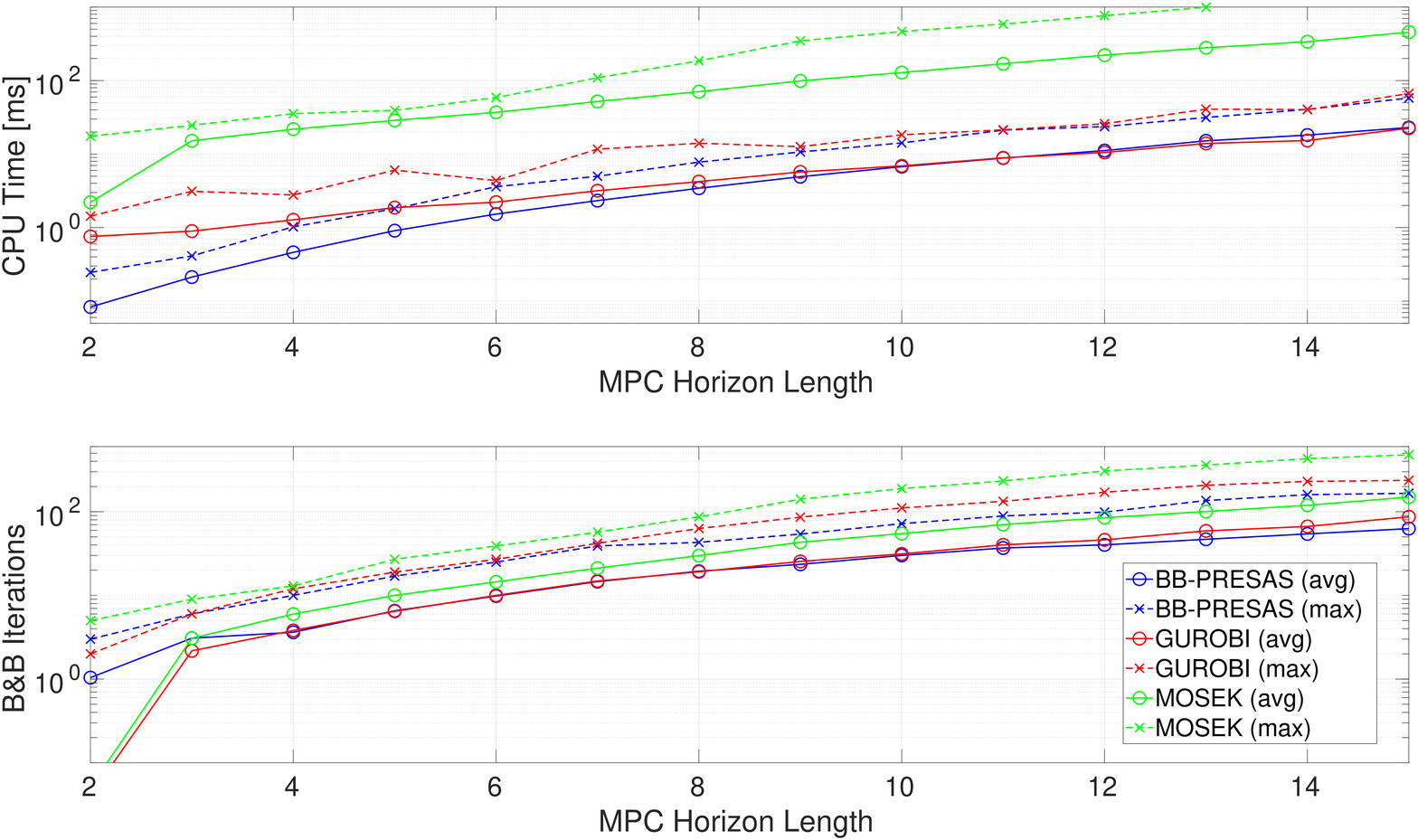}}}
\caption{Computational results for closed-loop mixed-integer MPC of the \texttt{bm99} example: \software{BB-PRESAS} versus \software{GUROBI} and \software{MOSEK} solvers for varying control horizon length $N$.}
\label{fig:BM99_time}
\end{figure}

The first case study is a hybrid MPC problem from~\cite{bemporad1999control}, with the default settings as in \texttt{bm99sim.m}, which is a part of the Hybrid Toolbox for \matlab. This demo example has been used also more recently for numerical comparisons in~\cite{Bemporad2018}. The system is modeled using the HYSDEL toolbox~\cite{torrisi2004hysdel} to obtain the mixed logical dynamical~(MLD) system formulation. 
%Figure~\ref{fig:BM99_traj} display the state trajectory and mode switching of the system. 
Figure~\ref{fig:BM99_time} illustrates the average and worst-case CPU times taken by our algorithm, \software{GUROBI} and \software{MOSEK} for a range of control horizon lengths $N$. 

Table~\ref{tab:hybrid_comp} presents a detailed comparison for this test example, including additional timing results for the \software{MI-NNLS} solver that are taken directly from~\cite{Bemporad2018}. The latter computational results can serve only as a reference since they have been obtained on a different computer, with respect to the one used here with a 2.80~GHz Intel Xeon E3-1505M v5 processor and 32~GB of RAM.
An important feature of our method is that its worst-case computation time is often rather close to the average performance in closed-loop MI-MPC simulations. This highlights the effectiveness of our tree propagation warm-starting procedure, such that consecutive branch-and-bound trees have approximately the same size.
In addition, it can be observed from Table~\ref{tab:hybrid_comp} that our proposed \software{BB-PRESAS} solver is either competitive with, or is a factor $2$ or $3$ times faster than \software{GUROBI}. The computational speedup is much larger when compared with other state-of-the-art tools such as \software{MOSEK}, our solver can be more than $10$ times faster in this particular MI-MPC test example.
It shall be noted that \software{GUROBI} is a heavily optimized and fairly large software, which is unlikely to be amenable for embedded microprocessors, due to its code size, memory requirements, and software library dependencies.

%\todo{add table with comparison against \software{GUROBI}, \software{OSQP}~\cite{stellato2018embedded} and \software{NNLS} solver~\cite{Bemporad2018}!}

%\begin{figure}[tpb] 
%\centering
%\includegraphics[trim={2.5in 0in 0.3in 0.4in},clip,scale=0.15]{BM99_traj}
%\caption{MPC state and input trajectory for BM99 system}
%\label{fig:BM99_traj}
%\end{figure}

\begin{table}[tpb]
	\caption{Timing results~(ms) per sampling step of hybrid MPC test problem for different horizon lengths $N$. Computation times for \software{MI-NNLS} solver are taken directly from~\cite{Bemporad2018}.}
	\label{tab:hybrid_comp}
	\centering
	\setlength{\tabcolsep}{0.75em}
	\begin{tabular}{l|c|c|c||c}
		\toprule
		N &{\software{BB-PRESAS}}&{\software{GUROBI}}&{\software{MOSEK}}&{\software{MI-NNLS}} \\
		& (mean/max) & (mean/max) & (mean/max) & (mean/max)\\
		\midrule
%		2 & 0.1/0.2 & 1.6/2.0 & 3.1/6.0 & 2.0/2.6 & -/- \\
%		3 & 0.4/0.6 & 2.7/3.0 & 6.3/11.5 & 2.5/4.8 & -/- \\
%		4 & 0.7/1.0 & 3.1/3.9 & 8.9/15.7 & 3.1/6.9 & -/- \\
%		5 & 1.1/2.3 & 3.9/4.8 & 10.8/15.7 & 3.9/13.0 & -/- \\
%		6 & 1.9/4.6 & 4.7/7.3 & 11.1/17.1 & 5.1/18.3 & -/- \\
%		7 & 3.0/8.9 & 5.6/9.5 & 17.5/80.4 & 6.4/30.2 & -/- \\
%		8 & 5.0/14.7 & 7.2/13.2 & 15.5/80.2 & 8.1/43.4 & -/- \\
%		9 & 7.5/20.7 & 8.8/15.3 & 22.8/110.6 & 11.1/69.8 & -/- \\
%		10 & 11.1/35.6 & 11.1/17.6 & 35.1/95.3 & 14.4/103.2 & -/- \\
%		11 & 15.7/48.3 & 13.0/23.9 & 37.3/102.5 & 20.6/179.1 & -/- \\
%		12 & 22.4/62.8 & 14.9/31.2 & 61.7/103.7 & 26.9/263.4 & -/- \\
%		13 & 29.5/75.0 & 17.5/36.6 & 47.3/119.5 & 35.5/384.9 & -/- \\
%		14 & 42.8/116.7 & 21.4/67.4 & 81.6/150.6 & 46.3/562.4 & -/- \\
%		15 & 50.4/133.0 & 25.9/109.8 & 89.9/181.1 & 61.7/766.9 & -/- \\
%		
		2 & 0.1/0.2 & 0.7/1.4 & 2.1/4.0 & 2.0/2.6 \\
		3 & 0.2/0.3 & 1.0/2.3 & 15.1/24.7 & 2.5/4.8 \\
		4 & 0.4/0.9 & 1.7/4.6 & 21.7/35.5 & 3.1/6.9 \\
		5 & 0.9/1.7 & 2.5/4.9 & 28.7/39.3 & 3.9/13.0 \\
		6 & 1.5/3.5 & 3.2/7.5 & 36.8/58.8 & 5.1/18.3 \\
		7 & 2.3/4.9 & 4.0/6.9 & 51.8/109.3 & 6.4/30.2 \\
		8 & 3.5/7.6 & 5.1/10.0 & 70.4/185.8 & 8.1/43.4 \\
		9 & 5.1/10.3 & 6.6/12.5 & 98.7/347.1 & 11.1/69.8 \\
		10 & 6.8/14.3 & 8.4/16.1 & 126.7/465.3 & 14.4/103.2 \\
		11 & 8.8/22.1 & 9.8/17.2 & 168.2/587.8 & 20.6/179.1 \\
		12 & 11.3/23.7 & 11.6/20.5 & 219.2/765.0 & 26.9/263.4 \\
		13 & 15.0/31.6 & 14.3/29.5 & 276.3/996.0 & 35.5/384.9 \\
		14 & 17.8/35.1 & 16.4/44.6 & 334.1/1241.9 & 46.3/562.4 \\
		15 & 21.0/41.6 & 21.9/71.6 & 450.8/1606.8 & 61.7/766.9 \\
		\bottomrule
	\end{tabular}
\end{table}

\subsection{Satellite Station Keeping with No-Go Zones} 

The second case study is motivated by a real-world application, namely, orbit control of a satellite in a circular low earth orbit, 400km from earth surface. The satellite propulsion system is composed of two on/off thrusters, one on each of the in-track faces of the satellite, with gimbals rotating along the vertical axis and subject to angle constraints~\cite{WalshWeissDiCairano16}. Thus, the propulsion system is controlled by two binary and two continuous control signals. The satellite dynamics are formulated by relative motion equations~(HCW) with respect to the target position along the orbit, and the cone constraints of the thrust forces are approximated as simplexes~\cite{WalshWeissDiCairano16}. Here, we consider a re-centering maneuver in which the satellite, previously drifting, is re-centered close to the target position along the orbit. Furthermore, the error coordinates from the target position are constrained in a station keeping window~($-300 \le X \le 300, -150 \le Y \le 150$).

%The second case study is motivated by the real-world application of controlling the in-plane motion of a satellite in order to keep track of a circular low earth orbit. We consider a satellite configuration with two gimbaled electric on-off thrusters, resulting naturally in a hybrid actuation system based on two binary and two continuous control input signals. The application and problem formulation is based on prior work as detailed in~\cite{Caverly2018}, even though we use a simplified problem formulation for the results in this paper. We additionally assume that the trajectory along such an orbit contains a no-go zone that needs to be avoided. This additional constraint is modeled using standard integer programming techniques (see~\cite{nemhauser1988integer}), resulting in three additional binary variables to implement the logical no-go zone constraints at each time point in the mixed-integer OCP. 

\begin{figure}[tpb] 
	\centerline{\hbox{
			\includegraphics[clip,width=0.55\textwidth]{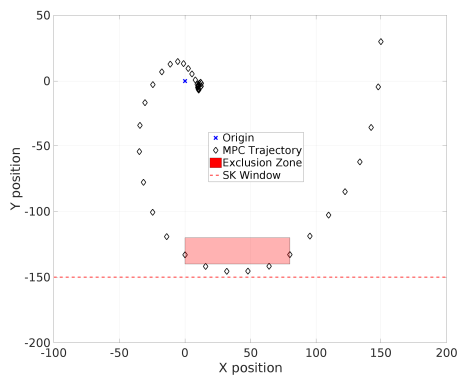}}}
\caption{MPC state evolution for satellite station keeping around the origin: rectangular no-go zone is depicted in red.}
\label{fig:zone_traj}
\end{figure}

\begin{figure}[tpb] 
	\centerline{\hbox{
			\includegraphics[width=0.55\textwidth]{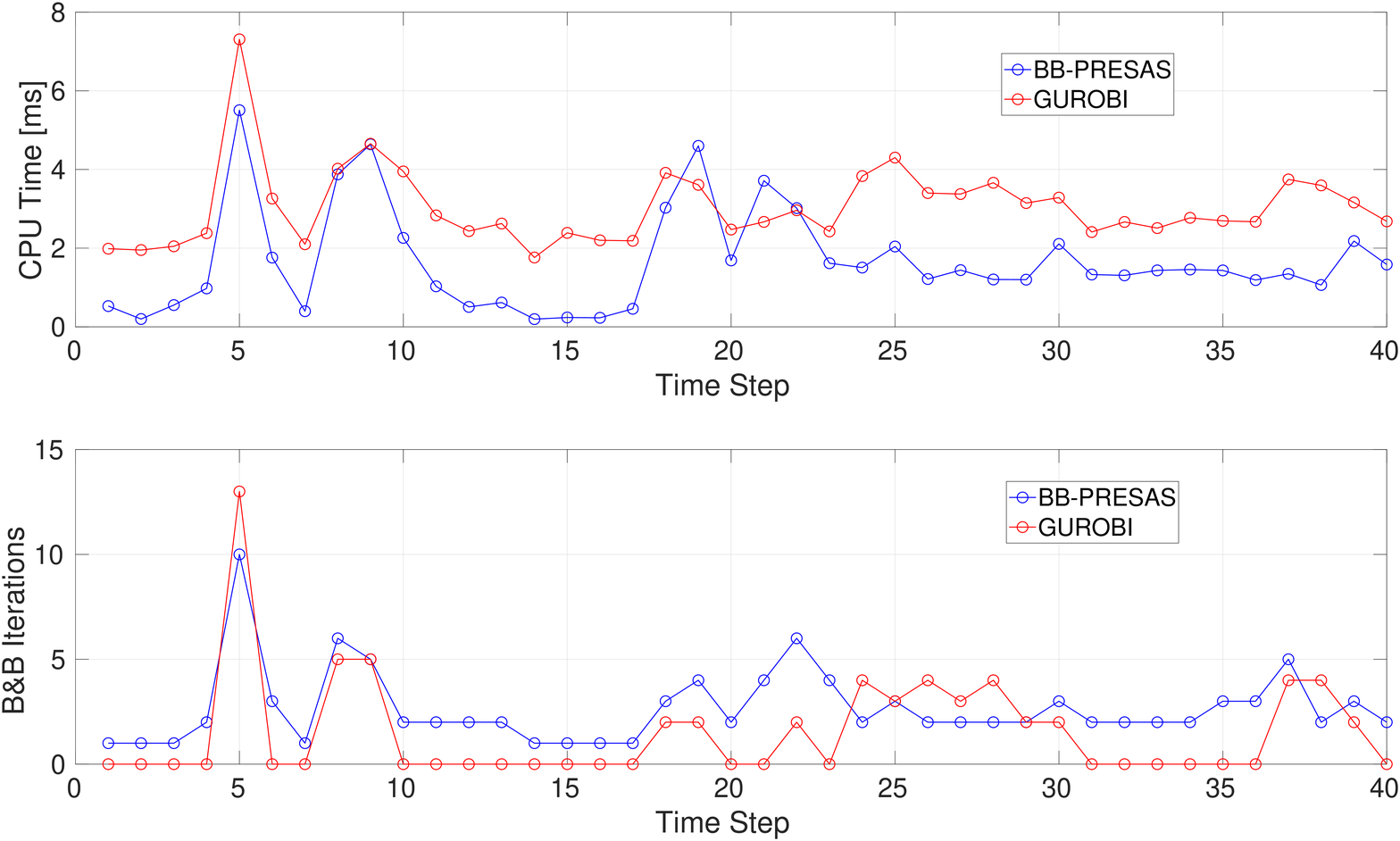}}}
\caption{Closed-loop results of mixed-integer MPC for satellite station keeping: \software{BB-PRESAS} versus \software{GUROBI} solver.}
\label{fig:zone_time}
\end{figure}

Thus, our problem is simplified from~\cite{WalshWeissDiCairano16}, by considering only the orbital dynamics in the orbital plane, i.e., ignoring the out-of-orbital-plane and attitude dynamics, and as a consequence using a simpler propulsion system with only two thrusters. To better highlight the potential of the MI-MPC method, we add an exclusion zone in the station keeping window, i.e., an area that must be avoided, which makes the allowed region of positions to be non-convex. This additional constraint is modeled using standard integer programming techniques (see, e.g.,~\cite{nemhauser1988integer}), resulting in three additional binary variables for each prediction step of the mixed-integer OCP to implement the logical exclusion zone constraints.
In Figure~\ref{fig:zone_traj}, we show the trajectory of the satellite in relative coordinates, where the origin is the desired satellite position along the orbit, for the simulation of the satellite controlled by the mixed-integer MPC.
%Hence when the system states reach the origin, it means the satellite is in the correct orbit. 
The depicted area in the figure corresponds to the station keeping window, in which the satellite should be kept, and the shaded area is the exclusion zone that must be avoided, at least pointwise in time. The computational timing results for this particular closed-loop MPC simulation can be found in Figure~\ref{fig:zone_time}. One can observe that our proposed algorithm has a very competitive runtime at every MPC time step, when compared to the commercial \software{GUROBI} solver. Most importantly, the \software{BB-PRESAS} algorithm appears to perform at least as good for this particular case study in terms of worst-case computation times. 

%\footnotetext{Computation times for \software{MI-NNLS} solver are taken directly from~\cite{Bemporad2018}.}

%\todo{Maybe a table with the breakdown in computation times of the C code?}

%%%%%%%%%%%%%%%%%%%%%%%%%%%%%%%%%%%%%%%%%%%%%%%%%%%%%%%%%%%%%%%%%%%%%%%%%%%%%%%%

%%%%%%%%%%%%%%%%%%%%%%%%%%%%%%%%%%%%%%%%%%%%%%%%%%%%%%%%%%%%%%%%%%%%%%%%%%%%%%%%
\section{CONCLUSIONS \& OUTLOOK}
\label{sec:concl}

In this paper, we proposed a branch-and-bound algorithm for mixed-integer MPC that exploits the optimal control problem structure to strengthen variable bounds, re-use pseudo-costs and warm-start the search tree at every MPC time step. More specifically, tailored domain propagation and tree propagation strategies have been presented. We showed preliminary results that illustrate the computational performance of our algorithm for two different MI-MPC case studies. 
%A more extensive set of simulations is part of ongoing work as to fully explore the design parameter space of the branch-and-bound method. 
%In addition, several other presolving features can be implemented efficiently in the MPC context, as well as alternative branching and node selection strategies. 
A compact, efficient, but self-contained C~code implementation of the proposed algorithm is under development to enable real-time embedded applications of hybrid MPC.

\bibliographystyle{IEEEtran}
\bibliography{IEEEabrv,references}

\end{document}